\documentclass[preprint,12pt]{elsarticle}

\newtheorem{theorem}{Theorem}[section]
\newtheorem{lemma}{Lemma}[section]
\newtheorem{corollary}{Corollary}[section]

\newtheorem{remark}{Remark}[section]
\newcommand{\ignore}[1]{}{}

\usepackage{amssymb}

\usepackage{color}
\usepackage{soul}
\usepackage[dvipsnames]{xcolor}

\usepackage[colorlinks=true, linkcolor=blue, citecolor=blue]{hyperref}

\def\1{{{\mbox{${\rm{1\negthinspace\negthinspace I}}$}}}}

\newcommand\beq{\begin{equation}}
\newcommand\eeq{\end{equation}}

\usepackage{ifthen}
\usepackage{xkeyval}
\usepackage{todonotes}
\begin{document}

\begin{frontmatter}

\title{Cram\'{e}r type moderate deviations for self-normalized $\psi$-mixing sequences }
\author{Xiequan Fan}
\address{Center for Applied Mathematics,
Tianjin University, Tianjin 300072,  China}

\begin{abstract}
Let $(\eta_i)_{i\geq1}$ be a  sequence of $\psi$-mixing random variables. Let $m=\lfloor n^\alpha \rfloor, 0< \alpha < 1, k=\lfloor n/(2m) \rfloor,$
and   $Y_j  = \sum_{i=1}^m \eta_{m(j-1)+i}, 1\leq j \leq k.$
Set $ S_k^o=\sum_{j=1}^{k } Y_j  $ and $[S^o]_k=\sum_{i=1}^{k } (Y_j )^2.$
We prove a Cram\'er type moderate deviation expansion for $\mathbb{P}(S_k^o/\sqrt{[ S^o]_k} \geq x)$ as $n\to \infty.$
Our result is similar to the recent work of Chen\textit{ et al.}\ [Self-normalized Cram\'{e}r-type moderate deviations under dependence.
Ann.\ Statist.\ 2016; \textbf{44}(4):  1593--1617]  where the authors established  Cram\'er type moderate deviation expansions
 for $\beta$-mixing sequences. Comparing to the result of Chen \textit{et al.}, our results hold  for
mixing coefficients with polynomial decaying rate and wider ranges  of validity.
\end{abstract}

\begin{keyword}  Cram\'{e}r moderate deviations, self-normalized processes, Studentized statistics, relative error, continued fraction expansions
\vspace{0.3cm}
\MSC Primary 62E20, 60F10, 60F15; Secondary  60G42
\end{keyword}

\end{frontmatter}




\section{Introduction}

The study of the relative errors for Gaussian approximations can be traced back to Cram\'{e}r (1938).
  Let  $(X_i)_{i\geq 1}$ be a sequence of independent and identically distributed (i.i.d.) centered real random variables satisfying the condition $\mathbb{E}\exp\{ c_{0}|X_{1}|\}<\infty$ for some constant $c_{0}>0.$
Denote $\sigma^2=\mathbb{E}X_{1}^2$ and $S_n=\sum_{i=1}^{n}X_{i}.$  Cram\'{e}r  established the following  asymptotic  moderate deviation expansion
on the tail probabilities of $S_n$:  For all $0\leq x =  o(n^{1/2}), $
\begin{equation}
\Bigg| \ln \frac {\mathbb{P}(S_n\geq x\sigma\sqrt{n})} {1-\Phi(x)} \Bigg|  = O(1)  \frac{(1+x)^3}{\sqrt{n}}  \ \ \mbox{as} \ \ n \rightarrow \infty,
\label{cramer1}
\end{equation}
where $\Phi(x)=\frac{1}{\sqrt{2\pi}}\int_{-\infty}^{x}\exp\{-t^2/2\}dt$ is the standard normal distribution.
In particular, inequality (\ref{cramer1}) implies that
\begin{equation}\label{tjk20}
\frac{\mathbb{P}(S_n\geq x\sigma\sqrt{n})}{1-\Phi \left( x\right)} = 1+ o(1)
\end{equation}
uniformly for $0\leq x =o( n^{ 1/6 } ) .$ Following the seminal work of Cram\'{e}r, various moderate deviation expansions for standardized  sums  have been obtained by many authors  (see, for instance, Petrov, 1954, 1975;   Linnik, 1961;   Saulis and Statulevi\v{c}ius, 1978;  Fan, 2017).    See also Ra\v{c}kauskas (1990, 1995),
Grama (1997),    Grama and Haeusler (2000),   Fan \textit{et al.}\ (2013)  for martingales.

To establish  moderate deviation expansions type of (\ref{tjk20}) for $0\leq x =o( n^{ \alpha} ), \alpha>0$,   we should assume that the random variables have finite moments of any order, see Linnik (1961). The last assumption becomes too restrictive if we only have finite moments of order $2+\delta, \delta \in (0, 1]$.  Thought we still can obtain (\ref{tjk20}) via Berry-Esseen estimations, the
range cannot wider than $|x|=O(\sqrt{\ln n}), n\rightarrow \infty.$ To overcome this shortcoming, a new type Cram\'{e}r type  moderate deviations (CMD),
called  self-normalized CMD, has been developed by Shao (1999).
Instead of considering the  moderate deviations for standardized sums $S_n/\sqrt{n\sigma^2},$ Shao (1999)  considered the moderate deviations
for self-normalized sums  $W_n:=S_n/\sqrt{\sum_{i=1}^n X_i^2}.$
Comparing to the standardized counterpart, the range of Gaussian approximation for self-normalized CMD can be much wider range than
 its counterpart for standardized  sums under same finite moment conditions. Moreover, in practice one usually does not
 known the variance of $S_n.$ Even the latter can be estimated, it is still advisable to
 use self-normalized CMD for more user-friendly. Due to these significant advantages, the study of
CMD for self-normalized sums  attracts more and more attentions.
For more self-normalized CMD  for independent random variables,
we refer to, for instance, Jing,   Shao and Wang (2003)    and Liu, Shao and Wang (2013).
  We also refer to de la Pe\~{n}a,   Lai and Shao (2009) and  Shao and Wang (2013)    for recent developments in this area.
For closely related results, see also de la Pe\~{n}a (1999) and Bercu and Touati (2008) for exponential inequalities for self-normalized martingales.

Thought self-normalized CMD for independent random variables has been well study,
 there are only a few of results for
 weakly dependent random variables. One of the main results  in this field is due to  Chen \textit{et al.}\ (2016).
Let $(\eta_i )_{i\geq1}$ be a (may be non-stationary) sequence of random variables.
Set $ \alpha \in (0, 1).$ Let $m=\lfloor n^\alpha \rfloor$ and $k=\lfloor n/(2m) \rfloor$, where $ \lfloor  x \rfloor$ denote the integer  part of $x.$
Denote  $$Y_j    = \sum_{i=1}^m\eta_{2m(j-1)+i}, \ \  1\leq j \leq k.$$
So $Y_j = S_{2m(j-1),m}.$ Set $$ S_k^o  =\sum_{j=1}^{k } Y_j  \ \ \ \ \  \textrm{and}  \ \ \ \ \ \  [S^o ]_k=\sum_{j=1}^{k } (Y_j )^2.$$
Define the  interlacing  self-normalized  sums as follows
\begin{equation}\label{Isums}
W_n^o= S_k^o  /  \sqrt{[ S^o ]_k}.
\end{equation}
Let $\mathcal{F}_{j} $ and $\mathcal{F}_{j+k}^{\infty}$
be $\sigma$-fields generated respectively by $(\eta_i)_{i \leq j}$ and $ (\eta_i)_{i \geq j+k}.$
The sequence of random variables  $(\eta_i )_{i\geq1}$ is called $\beta$-\textit{mixing} if the mixing coefficient
\begin{eqnarray}\label{sbgkm02}
 \beta(n):=  \sup_j  \mathbb{E} \sup \{ \big|\mathbb{P}( A |\mathcal{F}_{j} ) -\mathbb{P}(A) \big| :\ A \in \mathcal{F}_{j+n}^{\infty}    \}\rightarrow 0
  \ \ \ \ \ \ \textrm{as}\    n\rightarrow \infty.
\end{eqnarray}
See Doukhan (1994). Write $$S_{k,m}=\sum_{i=k+1}^{k+m}\eta_i$$
  the block sums of $(\eta_i )_{i\geq1}$ for $k+1 \leq i \leq k+m.$
   Throughout the paper,   denote $c,$  probably supplied with some indices,
  a generic positive constant.
Assume that $(\eta_i )_{i\geq1}$ are centered, that is
\begin{eqnarray}
   \mathbb{E}\eta_i =0   \ \ \ \textrm{for all} \   i , \label{cosf3.10}
\end{eqnarray}
and that there exists a constant  $\nu \in  (0, 1] $ such that
\begin{eqnarray}
 \mathbb{E}|\eta_i|^{2+\nu} \leq  c_0^{2+\nu}  \label{cosdf3.11}
\end{eqnarray}
and
\begin{eqnarray}
 \mathbb{E}S_{k,m}^2 \geq c_1^2 m \ \ \ \  \textrm{for all} \  k\geq0, m\geq 1.\label{cosdf3.12}
\end{eqnarray}
By Theorem 4.1 of Shao and Yu (1996),  it known that  condition (\ref{cosdf3.11})  usually implies
 the following condition: there exists a constant  $\rho \in  (0, 1] $ such that
 \begin{eqnarray}\label{cosdf3.13}
     \mathbb{E}|S_{k,m}|^{2+\rho} \leq m^{1+\rho/2}c_2^{2+\rho} ,
 \end{eqnarray}
  provided that  that the mixing coefficient  has a
  polynomially decaying rate as $n\rightarrow \infty$.
  In  (\ref{cosdf3.13}), it is usually that  $\rho  <  \nu .$
Assume conditions  (\ref{cosf3.10})-(\ref{cosdf3.12}). Assume also that there exist positive constants $a_1,   a_2$ and $\tau$ such that
$$ \beta(n) \leq a_1 e^{-a_2 n^\tau} .  $$
Using $m$-dependent approximation, Chen \textit{et al.}\,(2016) proved that for any positive $\rho<  \nu,$
\begin{equation} \label{dunfo1}
\Bigg|\ln \frac{\mathbb{P}(W_n^o \geq x)}{1-\Phi \left( x\right)} \Bigg| \leq   c_{ \rho}    \Bigg(  \frac{ (1+ x)^{2+\rho }}{n^{(1-\alpha)\rho / 2  }  }   \Bigg)
\end{equation}
uniform for  $0\leq x  =o( \min\{ n^{(1-\alpha)/2 },  n^{\alpha \tau/2 } \}),$
where $c_{ \rho}$ depends only on $c_0, c_1, \rho,  a_1,   a_2$ and $\tau.$
In particular, it implies that
\begin{equation}\label{dunfo}
\frac{\mathbb{P}(W_n^o \geq x)}{1-\Phi \left( x\right)} = 1+ o(1)
\end{equation}
uniformly for $0\leq x =o( \min \{ n^{ (1-  \alpha)\rho/(4+2\rho)  } ,   n^{\alpha \tau/2 }  \} ) .$  Equality (\ref{dunfo}) implies that
the tail probabilities of $W_n^o$ can be uniformly approximated by the standard normal distribution for moderate $x$'s.
Such type of results play an important role in statistical inference of means, see Section 5 of Chen \textit{et al.}\ (2016) for   applications.
  Inspiring the proof of Chen \textit{et al.}\ (2016),  it is easy to see that
(\ref{dunfo1}) remains valid when the conditions  (\ref{cosf3.10})-(\ref{cosdf3.12}) are replaced by the slightly more general  conditions (\ref{cosf3.10}), (\ref{cosdf3.12}) and (\ref{cosdf3.13}).

In this paper, we are interested to extend the results of Chen \textit{et al.}\ (2016)  to $\psi$-mixing sequences, with conditions (\ref{cosf3.10}), (\ref{cosdf3.12}) and (\ref{cosdf3.13}). By Proposition 1 in Doukhan (1994),  it is  known that  $\psi$-mixing usually implies
 $\beta$-mixing.  However, the ranges of our results do  not depend on the mixing coefficients. Indeed, our ranges  of validity for (\ref{dunfo1})
 and  (\ref{dunfo})  are respectively $0\leq x  =o(n^{(1-\alpha)/2 })$  and $0\leq x =o(  n^{ (1-  \alpha)\rho/(4+2\rho)  } )$ as $n\rightarrow \infty$, which are the best possible even $(\eta_i )_{i\geq1}$ are independent.
Moreover, we show that (\ref{dunfo}) remains true if $\psi$-mixing coefficient $\psi(n)$ decays in a  polynomial decaying rate,
in contrast to  $\beta$-mixing sequences which does not share this property.
 For methodology, our approach is based on  martingale approximation and self-normalized Cram\'{e}r type moderate deviations  for martingales due to Fan \textit{et al.}\ (2018).

%

The paper is organized as follows. Our main results are stated and discussed in Section \ref{sec bsum}.
Applications and  simulation study are   given in Section \ref{sec3}.
Proofs of results are deferred to Section \ref{sec4}.

\section{Main results}\label{sec bsum}


Recall that  $\mathcal{F}_{j} $ and $\mathcal{F}_{j+k}^{\infty}$
be $\sigma$-fields generated respectively by $(\eta_i)_{i \leq j}$ and $ (\eta_i)_{i \geq j+k}.$
We say that $(\eta_i )_{i\geq1}$
is $\psi$-\textit{mixing} if the mixing coefficient
\begin{eqnarray}\label{sgkm02}
 \psi(n):=  \sup_j   \sup_A\{ \big|\mathbb{P}( A |\mathcal{F}_{j} ) -\mathbb{P}(A) \big|/\mathbb{P}(A):\ A \in \mathcal{F}_{j+n}^{\infty}    \}\rightarrow 0
  \ \ \ \ \ \ \textrm{as}\    n\rightarrow \infty.
\end{eqnarray}
See Doukhan (1994).
Our  main result is the following  self-normalized Cram\'{e}r type moderate deviations
for  $\psi$-mixing sequences.
\begin{theorem}\label{th3.3} Assume conditions (\ref{cosf3.10}), (\ref{cosdf3.12}) and (\ref{cosdf3.13}).
Set $ \alpha \in (0, 1).$ Let $m=\lfloor n^\alpha \rfloor$ and $k=\lfloor n/(2m) \rfloor$ be respectively the integers part of $n^\alpha$ and $n/(2m)$.
Denote
$$  \delta_n^2=    m \psi^2(m )   +   k\psi(m )   $$
and
$$  \gamma_n =  k^{1/2}  \psi^{1/2}(m ) + n\psi(m ) .  $$
Assume also that $\delta_n,  \gamma_n\rightarrow 0$ as $n\rightarrow \infty.$
\begin{itemize}
   \item[\emph{[i]}] If $\rho \in (0, 1)$, then
for all  $0\leq x  =o( n^{(1-\alpha)/2 } ),$
\begin{equation}\label{ggdfg01}
\Bigg|\ln \frac{\mathbb{P}(W_n^o \geq x)}{1-\Phi \left( x\right)} \Bigg| \leq   c_{ \rho}    \Bigg(  \frac{  x ^{2+\rho }}{n^{(1-\alpha)\rho / 2  }  }  + x^2 \delta_n^2 + (1+ x) \big( \frac{1}{n^{(1-\alpha)\rho(2-\rho)/8}(1+x^{\rho(2+\rho)/4})  } +  \gamma_n \big)\Bigg) ,
\end{equation}
where $c_{ \rho}$ depends only on $c_1, c_2 $ and $ \rho.$

   \item[\emph{[ii]}] If $\rho =1$, then
for all  $0\leq x  =o( n^{(1-\alpha)/2} ),$
\begin{equation}\label{ggdfg02}
\Bigg| \ln \frac{\mathbb{P}(W_n^o \geq x)}{1-\Phi \left( x\right)} \Bigg|\leq  c    \Bigg(  \frac{  x ^{3}}{n^{(1-\alpha)  / 2  }  }  + x^2 \delta_n^2 + (1+ x) \big( \frac{1}{n^{(1-\alpha)/8}(1+x^{3/4})  } + \frac{ \ln n}{n^{(1-\alpha)/2 }}  +  \gamma_n \big)\Bigg) ,
\end{equation}
where $c $ depends only on $c_1$ and $ c_2. $
 \end{itemize}
\end{theorem}

Notice that in the i.i.d.\,case, $W_n^o$ is a self-normalized sums of $k$ i.i.d.\,random variables, that is $(Y_i)_{1\leq i \leq k}.$
According to the classical result of Jing, Shao and Wang (2003), Cram\'{e}r type moderate deviations holds for $0\leq x =o(k^{1/2}).$
Since the last range is equivalent to the range $0\leq x  =o( n^{(1-\alpha)/2} )$,
the ranges of validity for  (\ref{ggdfg01})  and  (\ref{ggdfg02}) coincide with the case of i.i.d., and, therefore, it is the best possible.

The following MDP result is a   consequence of the last theorem.
\begin{corollary}\label{corollary02} Assume the conditions of Theorem \ref{th3.3}.
Let $a_n$ be any sequence of real numbers satisfying $a_n \rightarrow \infty$ and $a_n/ n^{(1- \alpha)/2} \rightarrow 0$
as $n\rightarrow \infty$. Then for each Borel set $B \subset \mathbb{R}$,
\begin{eqnarray*}
- \inf_{x \in B^o}\frac{x^2}{2} &\leq & \liminf_{n\rightarrow \infty}\frac{1}{a_n^2}\ln \mathbb{P}\bigg(\frac{1}{a_n} W_n^o  \in B \bigg) \nonumber\\
 &\leq& \limsup_{n\rightarrow \infty}\frac{1}{a_n^2} \ln \mathbb{P}\bigg(\frac{1}{a_n} W_n^o \in B \bigg) \leq  - \inf_{x \in \overline{B}}\frac{x^2}{2} \,,
\end{eqnarray*}
where $B^o$ and $\overline{B}$ denote the interior and the closure of $B$, respectively.
\end{corollary}

If $\psi(n) = O\big(   n^{-   (1+ \rho) /   \alpha   }   \big),$ then  $\delta_n^2 = o\big( n^{-(1-\alpha)\rho  / 2  }  \big)$
and $\gamma_n=o(n^{-(1-\alpha)\rho  / 2  } ) .$
The following corollary is  nonetheless worthy to state.
\begin{corollary}\label{co2}  Assume conditions (\ref{cosf3.10}), (\ref{cosdf3.12}) and (\ref{cosdf3.13}). Set $ \alpha \in (0, 1).$
Assume also that  $$\psi(n) = O\big( n^{-   (1+ \rho) /   \alpha   } \big)$$ as $ n\rightarrow \infty.$
\begin{itemize}
   \item[\emph{[i]}] If $\rho \in (0, 1)$, then
for all  $0\leq x  =o( n^{(1-\alpha)/2 } ),$
\begin{equation}\label{sf01}
 \Bigg| \ln \frac{\mathbb{P}(W_n^o \geq x)}{1-\Phi \left( x\right)} \Bigg| \leq  c_{ \rho}    \Bigg(  \frac{  x ^{2+\rho }}{n^{(1-\alpha)\rho / 2  }  }   +  \frac{1+x}{n^{(1-\alpha)\rho(2-\rho)/8}(1+x^{\rho(2+\rho)/4})  }  \Bigg),
\end{equation}
where $c_{ \rho}$ depends only on $c_1, c_2 $ and $ \rho.$
   \item[\emph{[ii]}] If $\rho =1$, then
for all  $0\leq x  =o( n^{(1-\alpha)/2} ),$
\begin{equation}\label{sf02}
 \Bigg| \ln \frac{\mathbb{P}(W_n^o \geq x)}{1-\Phi \left( x\right)} \Bigg| \leq   c     \Bigg(  \frac{  x ^{3}}{n^{(1-\alpha)  / 2  }  }   + (1+ x) \big( \frac{1}{n^{(1-\alpha)/8}(1+x^{3/4})  } + \frac{ \ln n}{n^{(1-\alpha)/2 }}  \big)\Bigg) ,
\end{equation}
where $c$ depends only on $c_1$ and $ c_2 .$
 \end{itemize}
 In particular, (\ref{sf01}) and (\ref{sf02}) together implies that for $\rho \in (0, 1],$
\begin{equation}\label{tphisn3}
\frac{\mathbb{P}(W_n^o >x)}{1-\Phi \left( x\right)}= 1+ o(1)
\end{equation}
uniformly for $0\leq x=o( n^{ (1-  \alpha)\rho/(4+2\rho)  }).$
\end{corollary}

Chen \textit{et al.}\ (2016) (see Section 3 therein) showed that if $\beta$-mixing coefficient $\beta(n)$ decays only polynomial slowly, then (\ref{tphisn3})
 is not valid at $x= (C \ln n)^{1/2}$ for sufficiently large constant $C.$ However,  Theorem  \ref{th3.3}
shows that   the range of validity of (\ref{tphisn3}) can be much wider when $\beta$-mixing is replaced by $\psi$-mixing.

Recall that in the i.i.d.\,case, $W_n^o$ is a self-normalized sums of $k$ i.i.d.\,random variables.
By Remark 2 of Shao (1999), the range of validity for (\ref{tphisn3}) is also the best possible.

\begin{remark}
  Notice that if $(\eta_i )_{i\geq1}$ satisfies conditions (\ref{cosf3.10}), (\ref{cosdf3.12}) and (\ref{cosdf3.13}), then
   $(-\eta_i )_{i\geq1}$ also  satisfies the same conditions.
   Thus the assertions in  Theorem  \ref{th3.3}  and
 Corollary  \ref{co2}  remain valid when $\displaystyle  \frac{\mathbb{P}(W_n^o \geq x)}{1-\Phi \left( x\right)}$ is replaced  by $\displaystyle  \frac{\mathbb{P}(W_n^o \leq -x)}{ \Phi \left( -x\right)}$.
\end{remark}

\section{Applications}\label{sec3}
\subsection{Application to simultaneous confidence intervals}
Consider the problem of constructing  simultaneous confidence intervals for the mean value $\mu$ of the random variables $(\zeta_i)_{i\geq 1}$.
Assume that $(\zeta_i-\mu)_{i\geq 1}$ satisfies  the  conditions  (\ref{cosf3.10}), (\ref{cosdf3.12}) and (\ref{cosdf3.13}).
Let
$$T_n= \frac{\sum_{j=1}^k ( Y_j - m\mu ) }{\sqrt{\sum_{j=1}^k (Y_j - \overline{Y}_j)^2}}, $$
where $m=\lfloor n^\alpha\rfloor, k=\lfloor n/(2m)\rfloor,$ $Y_j= \sum_{i=1}^m \zeta_{2m(j-1)+i}, \  1\leq j \leq k,$ and $\overline{Y}_j=k^{-1}\sum_{j=1}^k Y_j.$
\begin{corollary} \label{c0kls}
 Let $\delta_n \in (0, 1).$  Assume that
\begin{eqnarray}\label{keldet}
 \big| \ln \delta_n \big| =o \big(   n^{(1-\alpha)\rho/(2+ \rho)} \big) .
\end{eqnarray}
If $\psi(n) = O\big( n^{-   (1+ \rho) /   \alpha   } \big),\,n\rightarrow \infty,$
then
\begin{eqnarray*}
\frac{\sum_{j=1}^k Y_j }{k m}  \pm  \frac{\Phi^{-1}(1-\delta_n/2)  }{km } \sqrt{ \sum_{j=1}^k (Y_j - \overline{Y}_j)^2}
\end{eqnarray*}
is $1-\delta_n$ conservative simultaneous confidence intervals for $\mu.$
\end{corollary}
\emph{Proof.}
 It is known that for all $x\geq0,$
\[
\mathbb{P}\Big( T_n  \geq x \Big) = \mathbb{P}\left(  \frac{\sum_{j=1}^k ( Y_j - m\mu  ) }{\sqrt{\sum_{j=1}^k (Y_j -  m\mu )^2}}  \geq x \Big(  \frac{k}{k-1}  \Big)^{1/2} \Big(\frac{k}{k+x^2-1} \Big)^{1/2}  \right),
\]
see Chung (1946). The last  equality  and (\ref{tphisn3}) together implies that
\begin{equation} \label{tphisns4}
\frac{\mathbb{P}(T_n  \geq x)}{1-\Phi \left( x\right)}= 1+ o(1)
\end{equation}
uniformly for $0\leq x=o( n^{ (1-  \alpha)\rho/(4+2\rho)  }).$
 Clearly, the upper $(\delta_n/2)$th quartile of a standard normal distribution $\Phi^{-1}(1-\delta_n/2)$  satisfies
$$\Phi^{-1}(1-\delta_n/2) = O(\sqrt{| \ln \delta_n |} ), $$  which, by (\ref{keldet}), is of order $o(n^{(1-\alpha)\rho/(4+ 2\rho)}).$
Then applying the last equality to $T_n$, we complete the proof of  Corollary \ref{c0kls}.

\subsection{Application to continued fraction and   simulation study}
One of the well known example of $\psi$-mixing sequences is called continued fraction expansions of irrational numbers on $(0, 1).$
For an irrational number $x \in (0, 1),$ let
\begin{eqnarray*}
a_1(x)  = \lfloor 1/x\rfloor, \, a_{n+1} (x)= a_1( x\circ T^n  ) , \ \  \ \ \ \ n\geq 1,
\end{eqnarray*}
be the continued fraction expansion of $x,$  where $T$ is defined by $T(x) =1/x- \lfloor 1/x\rfloor$, that is the fractional part of $1/x.$
It is easy to see that
\begin{eqnarray*}
\displaystyle x=\frac{1}{\displaystyle  a_1(x) + \frac{1}{\displaystyle   a_2(x) +\frac{1}{\displaystyle   a_3(x) + \frac{1}{ \cdot \cdot \cdot } }  } }
\end{eqnarray*}
The sequence  $(a_n(x))_{n\geq 1}$ with respect to the uniform measure in $(0, 1)$  is $\psi$-mixing. Indeed, L\'{e}vy (1929)    proved that
\begin{eqnarray}
 \psi(n) =\sup_j   \sup_A\{ \big|\mathbb{P}( A |\mathcal{F}_{j} ) -\mathbb{P}(A) \big|/\mathbb{P}(A):\ A \in \mathcal{F}_{j+n}^{\infty}    \} \leq C e^{-\lambda n}
\end{eqnarray}
with positive absolute constants $C$ and $ \lambda$, where $\mathcal{F}_{1}^j $ and $\mathcal{F}_{j+ n}^{\infty}$
be $\sigma$-fields generated respectively by $(a_i(x))_{1 \leq i \leq j}$ and $ (a_i(x))_{i \geq  j+n}.$
Denote by
\begin{eqnarray*}
\mathbb{G}(E)  = \frac{1}{\ln 2} \int_{E} \frac{1}{1+x}dx ,
\end{eqnarray*}
the Gauss measure on the class of Borel subsets $\mathcal{B}$ of $(0, 1).$
 It is known that (cf.\ Billingsley (1965)) $T$ is
an ergodic transformation preserving the Gauss measure and thus  $(a_n(x))_{n\geq 1}$ is a stationary ergodic sequence  with respect to the probability space
$((0, 1), \mathcal{B}, \mathbb{G})$. Clearly,
the set $\{a_1=k \}$ is the interval $(1/(k+1), 1/k]$ and thus
\begin{eqnarray}
\mathbb{G}(\{a_1=k\})= \frac{1}{\ln 2} \int_{1/(k+1)}^{1/k}  \frac{1}{1+x}dx= \frac{1}{\ln 2} \ln\!\Big(1+\frac{1}{k(k+2)} \Big)   . \nonumber
\end{eqnarray}
Hence, by the ergodic theorem we have for any function $F: \mathbb{N}\rightarrow \mathbb{R},$ it holds
\begin{equation}\label{lksed}
 \lim_{N\rightarrow \infty}\frac{1}{N}  \sum_{k=1}^N F(a_k(x))=    \frac{1}{\ln 2} \sum_{j=1}^{\infty} F(j) \ln\!\Big(1+\frac{1}{j(j+2)} \Big) \ \ \ \textrm{a.e.}
\end{equation}
whenever the series on the right hand side converges absolutely.
Recently,  Bazarova, Berkes and   Horv\'{a}th (2016)   gave a central limit theorem for $(a_n(x))_{n\geq 1}$.
Next, we give a self-normalized  Cram\'{e}r type moderate deviations.

Letting $\mathbb{E}$ denote expectation with respect to $\mathbb{G}$, by (\ref{lksed}), we have $\mathbb{E} a_1(x)= \infty$ and $\mathbb{E} (a_1(x))^\alpha <  \infty$ for any $\alpha \in (0, 1).$
Consider the self-normalized moderate deviation for the random variables $(\zeta_i )_{i\geq 1}$,  where $\zeta_i=   \sqrt[3]{a_i(x)} $ for any $i$.
Then $\mathbb{E}(\zeta_1)^{2+\rho} < \infty$ for any $\rho \in (0 , 1)$ and
\begin{equation}
\mu:=\mathbb{E}\zeta_i =       \frac{1}{\ln 2} \sum_{j=1}^{\infty}  j^{1/3} \ln\!\Big(1+\frac{1}{j(j+2)} \Big).
\end{equation}
Let
$$W_n^o= \frac{\sum_{j=1}^k ( Y_j - m\mu  ) }{\sqrt{\sum_{j=1}^k (Y_j - m\mu )^2}}, $$
where $m=\lfloor n^\alpha\rfloor, k=\lfloor n/(2m)\rfloor,$ $Y_j= \sum_{i=1}^m \zeta_{2m(j-1)+i}, \  1\leq j \leq k.$
By (\ref{tphisn3}), we have
the following result.
\begin{corollary}\label{c041}
Set $ \alpha \in (0, 1).$ Then for any $\rho \in (0, 1),$
\begin{equation}
\frac{\mathbb{P}(W_n^o \geq t)}{1-\Phi \left( t\right)}= 1+ o(1)
\end{equation}
uniformly for $0\leq t=o( n^{ (1-  \alpha)\rho/(4+2\rho)  }).$
\end{corollary}

Next, we give a simulation study for the last corollary.
We let $n=30$, $m= 1, 2, 3,4$ and consider 13 levels of $t: t=0, .1, .2, ..., 1.0, 1.2, 1.4.$
Let $x$ be the discrete uniform distribution  random variable, with possible values  $ \pi/10000, 2\pi/10000,....., 3182\pi/10000 .$  Since $\pi$ is an irrational number,  $x$ are   irrational numbers.  In $W_n^o,$
we take
\begin{equation}
\mu =   \frac{1}{\ln 2} \sum_{j=1}^{300}  j^{1/3} \ln\!\Big(1+\frac{1}{j(j+2)} \Big) .
\end{equation}
Then $\mathbb{P}(W_n^o \geq t)\approx \#(W_n^o: W_n^o \geq t) / 3182. $ The following table shows the simulate rations $\frac{\mathbb{P}(W_n^o \geq t)}{1-\Phi \left( t\right)}.$ From the table, we see that   the  interlacing  self-normalized  sums (that is $m=2,3, 4$)    has a better performance than
  self-normalized  sums (that is $m=1$) when $x$ close to $0$. When $x$ moves away from $0,$ the reverse is true.
\begin{center}
\begin{tabular}{ c l l l l l l l l l l  l l  l  }
\hline
$m$ &   $t=0$ &   $ .1$ &  $.2$ &   $.3$ &   $.4$ &    $.5$ &    $.6$ &    $.7$ &   $.8$&   $.9$&  $1.0$&   $1.2$  &   $1.4$ \\
\hline
  $1$ & $1.11$ & $1.13$ & $1.15$&  $1.16$&  $1.17$ &  $1.16$ &  $1.11$ & $1.08$& $1.03$&  $0.96$   & $0.90$   &  $0.75$    &  $0.53$   \\
  $2$ & $1.01$ & $ 1.02$ & $1.02$&  $1.02$&  $1.02$ &  $1.01$ &  $1.00$ & $0.99$& $0.94$&  $0.88$ & $0.78$   &   $0.57$     &  $0.42$            \\
  $3$ & $1.00 $ & $1.03$ & $1.04$&  $1.07$&  $1.06$ &  $1.06$ &  $1.04$ & $1.01$& $0.98$&  $0.92$  & $0.85$   &   $0.67$    &  $0.48$    \\
  $4$ & $1.01 $ & $1.00$ & $0.99$&  $0.96$&  $0.94$ &  $0.89$ &  $0.82$ & $0.74$& $0.67$&  $0.56$  & $0.46$   &   $0.29$    &  $0.13$    \\
\hline
\end{tabular}
\end{center}

\section{Proofs} \label{sec4}
To shorten notations, for two real positive sequences $(a_n)_{i\geq1}$ and $(b_n)_{i\geq1},$ write $a_n \preceq b_n $ if there exists a positive constant $C$
such that $ a_n  \leq C  b_n $ holds for all large $n,$ $a_n \succeq b_n $ if   $b_n \preceq a_n,$ and $a_n \asymp b_n $ if $  a_n \preceq b_n $ and $b_n \preceq a_n.$

We only give a proof for the case where $\rho \in (0, 1).$ For the case where $\rho=1,$ the proof is similar.
\subsection{Preliminary lemmas }

Let $(X_i,\mathcal{F}_i)_{i=0,...,n} $  be a sequence of   martingale differences defined on a
 probability space $(\Omega ,\mathcal{F},\mathbb{P})$. Set
\begin{equation}
S_{0}=0,\ \ \ \ \ S_k=\sum_{i=1}^k X_i,\quad k=1,...,n.  \label{xk}
\end{equation}
Then $(S_k,\mathcal{F}_k)_{k=0,...,n}$ is a martingale. Denote $B_n^2=\sum_{i=1}^n\mathbb{E}X_i^2$ the variance of $S_n$.
We assume the following conditions:
\begin{description}
\item[(A1)]  There exists   $ \varsigma_n \in [0, \frac14] $ such that
$$ \Big| \sum_{i=1}^n\mathbb{E}[ X_{i}^2  | \mathcal{F}_{i-1}] - B_n^2 \Big| \leq  \varsigma_n^2 B_n^2 ;$$
\item[(A2)]  There exist  $\rho \in (0, 1]$ and   $ \tau_n \in (0, \frac14]  $ such that
\[
\mathbb{E}[|X_{i}| ^{2+\rho}  | \mathcal{F}_{i-1}]   \leq  (\tau_{n }  B_n)^{\rho}\, \mathbb{E}[ X_{i}  ^{2}  | \mathcal{F}_{i-1}], \ \ \ 1\leq i \leq n.
\]
\end{description}
In practice, we usually have $ \varsigma_n, \tau_{n}\rightarrow 0$ as $n \rightarrow \infty$. In the case of sums of
i.i.d.\,random variables with finite $(2+\rho)$th moments, then it holds $B_n\asymp \sqrt{n}$, and thus conditions (A1)  and (A2) are satisfied with $ \varsigma_n=0$ and $\tau_{n }=O(1/\sqrt{n} )$ as $n\rightarrow \infty.$

Define the self-normalized martingales
\begin{equation}
W_n = \frac{ S_n}{\sqrt{\sum_{i=1}^n X_i^2}} .
\end{equation}
The proof of Theorem \ref{th3.3} is based on the  following technical lemma due to Fan \textit{et al.}\ (2018) (see Corollary 2.3 therein), which gives a Cram\'er type moderate deviation expansion for self-normalized martingales.
\begin{lemma}\label{lemmpr1}
Assume  conditions (A1) and (A2).  Denote
 \begin{equation}  \label{alpphahat001}
 \widehat{\tau}_n(x, \rho) =  \frac{ \tau_{n }^{ \rho(2-\rho)/4 } }{ 1+ x  ^{  \rho(2+\rho)/4 }}.
 \end{equation}
\begin{description}
\item[\textbf{[i]}] If $\rho \in (0, 1)$, then for $0\leq x =o(\tau_n^{-1})$,
\begin{eqnarray*}
\Bigg|  \ln \frac{\mathbb{P}(W_n \geq x)}{1-\Phi \left( x\right)}  \Bigg| \leq  c_{\rho} \Bigg( x^{2+\rho}  \tau_n^\rho+ x^2 \varsigma_n^2 +(1+x)\Big(
    \varsigma_n +  \widehat{\tau}_n(x, \rho)\Big) \Bigg) ,
\end{eqnarray*}
where $c_\rho$ depends only on $\rho.$
\item[\textbf{[ii]}] If $\rho =1$, then for $0\leq x =o(\gamma_n^{-1})$,
\begin{eqnarray*}
 \Bigg| \ln \frac{\mathbb{P}(W_n \geq x)}{1-\Phi \left( x\right)} \Bigg|   \leq c  \Bigg( x^{3}  \tau_n + x^2 \varsigma_n^2+(1+x)\Big(  \varsigma_n+\tau_n |\ln \tau_n |+ \widehat{\tau}_n(x, 1)  \Big) \Bigg)  ,
\end{eqnarray*}
where $c$ is a constant.
\end{description}
\end{lemma}

The following lemma  is useful in the proof of Theorem \ref{th3.3}, see Theorem 2.2 of Fan \textit{et al.}\ (2017).
Denote   $x^+=\max\{x, 0\}$ and $x^-=(-x)^+ $  the positive and negative parts of $x,$ respectively.
\begin{lemma}
\label{th24} Assume that $\mathbb{E} |X _i|^{\beta} < \infty$ for a constant $\beta \in (1, 2]$ and  all $i\in [1,n]$. Write
\[
\textrm{G}^0_k(\beta)  =\sum_{i=1}^k \Big(  \mathbb{E}\big[ (X_i^-)^\beta  |\mathcal{F}_{i-1} \big] + (X_i^+)^\beta \Big),\ \ \ \ k\in [1,n].
\]
Then  for all $x, v>0$,
\begin{eqnarray}
 \mathbb{P}\left( S_k  \geq x\ \mbox{and}\ \textrm{G}^0_k(\beta) \leq v^\beta \ \mbox{for some}\ k\in[1,n]\right)
 \leq   \exp\left\{-C(\beta) \left(\frac{x}{v} \right)^\frac{\beta}{\beta -1} \right\}, \label{f71}
\end{eqnarray}
where
\begin{eqnarray}\label{cbeta}
C(\beta)= \beta ^{\frac1{1-\beta}} \left(1- \beta^{-1} \right).
\end{eqnarray}
\end{lemma}

In the proof of Theorem \ref{th3.3}, we also make use of the following lemma which can be found in  Theorem 3 of Doukhan (1994).
\begin{lemma}\label{HHA6}
Suppose that $X$ and $Y$ are random variables which are $\mathcal{F}_{j+n}^{\infty}$- and $\mathcal{F}_{j}$-measurable, respectively, and that $\mathbb{E}|X| < \infty,$ $\mathbb{E}|Y|  < \infty$.
Then
$$\Big|\mathbb{E}  X  Y  - \mathbb{E}  X\mathbb{E} Y \Big|\leq  \psi(n)  \,  \mathbb{E}|X|  \,  \mathbb{E}|Y| .$$
Moreover, since $\mathbb{E}|X| \leq (\mathbb{E}|X|^2)^{1/2},$  it holds
$$\Big|\mathbb{E}  X  Y  - \mathbb{E}  X\mathbb{E} Y \Big|\leq  \psi(n)  \,  (\mathbb{E} X^2)^{1/2}  \, ( \mathbb{E} Y^2)^{1/2} $$
provided that $\mathbb{E} X^2 < \infty$ and $\mathbb{E} Y^2  < \infty$.
\end{lemma}

\subsection{Proof of Theorem \ref{th3.3}}
Denote by $\mathcal{F}_{l} = \sigma \{\eta_{i}, 1\leq i \leq  2ml-m  \}.   $ Then $Y_j$ is $\mathcal{F}_{j}$-measurable.
Since $\mathbb{E} \eta_{i } =0$ for all $i$, by the definition of mixing coefficient (\ref{sgkm02}) and condition (\ref{cosdf3.13}),  it is easy to see that for $1\leq j \leq k,$
\begin{eqnarray}
 \Big|\mathbb{E}[ Y_j | \mathcal{F}_{ j-1 } ]\Big|&=& \Big|\sum_{i=1 }^{ m  } \Big( \mathbb{E}[\eta_{2m(j-1)+i } | \mathcal{F}_{j-1} ]- \mathbb{E} \eta_{2m(j-1)+i }  \Big)\Big| \nonumber \\
  &\leq& \sum_{i=1 }^{ m  } \psi(m) \mathbb{E} |\eta_{2m(j-1)+i }|   \nonumber  \\
  &\leq& \sum_{i=1 }^{ m  } \psi(m) (\mathbb{E} |\eta_{2m(j-1)+i }|^{2+\rho})^{1/(2+\rho)}   \nonumber  \\
   &\leq& m\psi(m) c_2,   \label{ineq01}
\end{eqnarray}
where  $c_2$ is given by (\ref{cosdf3.13}).
Thus
\begin{eqnarray*}
\Big|\sum_{j=1}^k \mathbb{E}[ Y_j | \mathcal{F}_{ j-1 } ]  \Big| \leq  k m \psi(m) c_2 \leq n \psi(m ) c_2  .
\end{eqnarray*}
By condition (\ref{cosdf3.13}) and the inequality
\begin{eqnarray*}
(x+y)^p \leq 2^{p-1}(x^p+y^p) \ \ \ \ \textrm{for}\  x, y\geq 0\ \textrm{and}\ p\geq 1,
\end{eqnarray*} we have
\begin{eqnarray}
 \mathbb{E}[  |Y_j-\mathbb{E}[ Y_j | \mathcal{F}_{ j-1 } ] |^{2+\rho} | \mathcal{F}_{ j-1 } ]&\leq& 2^{1+\rho}\mathbb{E}[  |Y_j|^{2+\rho} + |\mathbb{E}[ Y_j | \mathcal{F}_{ j-1 } ] |^{2+\rho} | \mathcal{F}_{ j-1 } ] \nonumber  \\
 &\leq& 2^{2+\rho}\mathbb{E}[  |Y_j|^{2+\rho}   | \mathcal{F}_{ j-1 } ] \nonumber \\
 &\leq& 2^{2+\rho} (1+ \psi(m )) \mathbb{E}  |Y_j|^{2+\rho} \nonumber  \\
  &\leq& 2^{2+\rho} (1+ \psi(m )) m^{1+\rho/2}c_2^{2+\rho}. \label{ine11.1}
\end{eqnarray}
The last inequality  implies that
\begin{eqnarray}
 \mathbb{E}[  |Y_j-\mathbb{E}[ Y_j | \mathcal{F}_{ j-1 } ] |^{2} | \mathcal{F}_{ j-1 } ]&\leq& ( \mathbb{E}[  |Y_j-\mathbb{E}[ Y_j | \mathcal{F}_{ j-1 } ] |^{2+\rho} | \mathcal{F}_{ j-1 } ] )^{2/(2+\rho)} \nonumber  \\
  &\leq& 2^{2} (1+ \psi(m ))^{2/(2+\rho)} m c_2^{2}\nonumber  \\
  &\leq& 2^{2} (1+ \psi(m ))  m c_2^{2}. \label{ine11.11}
\end{eqnarray}
Similarly,  by (\ref{cosdf3.13}) and the assumption  $ \delta_n \rightarrow 0$ as $n\rightarrow \infty$, it holds
\begin{eqnarray}
  \mathbb{E}[  |Y_j-\mathbb{E}[ Y_j | \mathcal{F}_{ j-1 } ] |^2 | \mathcal{F}_{ j-1 } ]
  &=&   \mathbb{E}[  Y_j^2 | \mathcal{F}_{ j-1 } ]    -  | \mathbb{E}[ Y_j | \mathcal{F}_{ j-1 } ] |^2  \nonumber \\
 &\geq& (1- \psi(m )) \mathbb{E}  Y_j^2   - (m\psi(m) c_2)^2 \nonumber \\
&\succeq&   \frac12 (1 - \psi(m )) m  c_1^2 .\label{ine11.2}
\end{eqnarray}
Combining (\ref{ine11.1})-(\ref{ine11.2}), we deduce that
$$ \mathbb{E}[  |Y_j-\mathbb{E}[ Y_j | \mathcal{F}_{ j-1 } ] |^{2+\rho} | \mathcal{F}_{ j-1 } ]\preceq m^{\rho/2}  \mathbb{E}[  |Y_j-\mathbb{E}[ Y_j | \mathcal{F}_{ j-1 } ] |^2 | \mathcal{F}_{ j-1 } ] ,  $$
$$ \sum_{j=1}^k \mathbb{E}[  |Y_j-\mathbb{E}[ Y_j | \mathcal{F}_{ j-1 } ] |^2 | \mathcal{F}_{ j-1 } ] \asymp  n $$
and, by Lemma \ref{HHA6} and (\ref{ineq01}),
\begin{eqnarray*}
&&\Big|\sum_{j=1}^k \mathbb{E}[  |Y_j-\mathbb{E}[ Y_j | \mathcal{F}_{ j-1 } ] |^2 | \mathcal{F}_{ j-1 } ] - \mathbb{E}  S_n^2 \Big| \\
 &&\leq\Big|\sum_{j=1}^k \mathbb{E}[  |Y_j-\mathbb{E}[ Y_j | \mathcal{F}_{ j-1 } ] |^2 | \mathcal{F}_{ j-1 } ] -\sum_{j=1}^k\mathbb{E}  Y_j^2 \Big| +\Big|  \mathbb{E}  S_n^2-\sum_{j=1}^k\mathbb{E}  Y_j^2 \Big|\\
  &&\leq  \sum_{j=1}^k \Big| \mathbb{E}[  Y_j^2 | \mathcal{F}_{ j-1 } ] -\mathbb{E}  Y_j^2 \Big|     + \sum_{j=1}^k \Big| \mathbb{E}[ Y_j | \mathcal{F}_{ j-1 } ] \Big|^2 + 2\sum_{j= 1}^k \sum_{l=1 }^{j-1}\Big|\mathbb{E}  Y_j  Y_l   \Big|\\
 &&\leq k\psi(m )\mathbb{E}Y_j^2 + k ( m\psi(m) c_2 )^2+ 2 \psi  (m ) \sum_{j= 1}^k \sum_{l=1 }^{j-1}  \mathbb{E}|Y_j|\,  \mathbb{E}|Y_l|  \\
  &&\leq 2 n \psi(m )c_2^2+n m \psi^2(m ) c_2^2 + 2 \psi  (m ) \sum_{j= 1}^k \sum_{l=1 }^{j-1}  \sqrt{\mathbb{E} Y_j^2}\, \sqrt{ \mathbb{E} Y_l^2}  \\
 &&\leq 2 n \psi(m )c_2^2+n m \psi^2(m ) c_2^2 + 2n k \psi (m ) c_2^2 \\
 && \leq  n m \psi^2(m ) c_2^2 + 4n k  \psi(m )c_2^2  .
\end{eqnarray*}
Denote by
$$\delta_n^2=   m \psi^2(m )   +   k\psi(m )  .$$ Taking $X_i=Y_j-\mathbb{E}[ Y_j | \mathcal{F}_{ j-1 } ], $ we find that  condition (A1) and (A2) is satisfied with $B_n^2=\mathbb{E}  S_n^2 \asymp  n,$ $ \varsigma_n\asymp \delta_n$ and $\tau_n \asymp \sqrt{m/n} \asymp    n^{-(1-\alpha)/2} .$
Applying Lemma  \ref{lemmpr1}  to
$$W_n:= \frac{  \sum_{j=1}^k (Y_j-\mathbb{E}[ Y_j | \mathcal{F}_{ j-1 } ])  }{ \sqrt{ \sum_{j=1}^k (Y_j-\mathbb{E}[ Y_j | \mathcal{F}_{ j-1 } ])^2}\ },$$
we have for all $0\leq x =o(n^{(1-\alpha)/2})  ,$
\begin{equation}\label{ineq413}
\Bigg| \ln \frac{\mathbb{P}(W_n \geq x)}{1-\Phi \left( x\right)} \Bigg| \leq c_{ \rho} \,  \Bigg(  \frac{  x ^{2+\rho }}{n^{(1-\alpha)\rho / 2  }  }  + x^2 \delta_n^2 + (1+ x) \Big( \frac{1}{n^{(1-\alpha)\rho(2-\rho)/8}(1+x^{\rho(2+\rho)/4})  } +  \delta_n \Big)\Bigg) .
\end{equation}
Notice that, by Cauchy-Schwarz's inequality,
\begin{eqnarray}
 \Big| \sum_{j=1}^k \big(  Y_j-\mathbb{E}[ Y_j | \mathcal{F}_{ j-1 } ]  \big)^2- \sum_{j=1}^k Y_j^2  \Big|
 &\leq&  2\sum_{j=1}^k |Y_j \mathbb{E}[ Y_j | \mathcal{F}_{ j-1 } ] | +  \sum_{j=1}^k (\mathbb{E}[ Y_j | \mathcal{F}_{ j-1 } ])^2 \nonumber  \\
&\leq& 2  m \psi(m) c_2  \sum_{j=1}^k |Y_j   | +  \sum_{j=1}^k \big (m \psi(m) c_2 \big)^2 \nonumber  \\
&\leq& 2 k^{1/2}    m \psi(m) c_2 \Big(\sum_{j=1}^k  Y_j^2 \Big)^{1/2}     +  k m^2 \psi^2(m) c_2^2.\ \ \ \ \ \ \   \label{sud01}
\end{eqnarray}
When $\sum_{j=1}^k  Y_j^2 \geq 1/4,$ both sides of the last inequality divided by $\sum_{j=1}^k  Y_j^2$,   we get
\begin{eqnarray*}
 \Bigg| \frac{\sum_{j=1}^k \big(  Y_j-\mathbb{E}[ Y_j | \mathcal{F}_{ j-1 } ]  \big)^2}{\sum_{j=1}^k Y_j^2}   -1  \Bigg|  \leq  4 k^{1/2}   m\psi(m ) c_2   + 4k m^2 \psi^2(m) c_2^2.
\end{eqnarray*}
By assumption $\gamma_n \rightarrow 0,$  we have $ k^{1/2}   m\psi(m )  \rightarrow 0$  as $n\rightarrow \infty$.
By Cauchy-Schwarz's inequality, we have $ \sum_{j=1}^k  |Y_j| \leq k^{1/2} \sqrt{ \sum_{j=1}^k   Y_j^2}  .$
Hence, when $\sum_{j=1}^k  Y_j^2 \geq 1/4,$  it holds
\begin{eqnarray*}
\Big| W_n - W_n^o  \Big| &=&  \Big| W_n \sqrt{\Sigma_{j=1}^k  Y_j^2} - \sum_{j=1}^k  Y_j  \Big| \frac{1}{ \sqrt{\sum_{j=1}^k  Y_j^2}  }  \\
&\preceq&  \sum_{j=1}^k\Big| \mathbb{E}[ Y_j | \mathcal{F}_{ j-1 } ]  \Big| + \frac{  \sum_{j=1}^k  |Y_j|  }{ \sqrt{\sum_{j=1}^k  Y_j^2}  } \Big(4 k^{1/2}   m\psi(m ) c_2   + 4k m^2 \psi^2(m) c_2^2 \Big)  \\
 &\leq&   n \psi(m ) c_2  +     k^{1/2} \Big (4 k^{1/2}   m\psi(m ) c_2   + 4k m^2 \psi^2(m) c_2^2 \Big) \\
 & =&9  n\psi(m ) c_2 .
\end{eqnarray*}
Hence,  when $\sum_{j=1}^k  Y_j^2 \geq 1/4,$  we have $$\Big| W_n - W_n^o  \Big| \leq  C \varepsilon_n,$$
where $C$ is a positive constant and $$\varepsilon_n =    n\psi(m ) c_1.$$
Notice that for $x\geq 0$ and $ |\varepsilon_n| =O(1) $,
\[
\frac{1-\Phi \left( x + \varepsilon_n \right)}{1-\Phi \left( x\right) }=   \exp\Big\{ O(1)(1+x) |\varepsilon_n| \Big\}.
\]
Without loss of generality, we may assume that $\sum_{j=1}^k \mathbb{E}Y_j^2=n;$ otherwise, we may consider   $(\eta_i/  \sqrt{\sum_{j=1}^k \mathbb{E}Y_j^2/n}  )_{1\leq i \leq n}$ instead of $(\eta_i)_{1\leq i \leq n}$.
Then it follows that
\begin{eqnarray}
 \mathbb{P}\Bigg(   \sum_{j=1}^k  Y_j^2 < \frac1 4 \Bigg)
&\leq& \mathbb{P}\Bigg(   \sum_{j=1}^k    (  Y_j^2 -\mathbb{E}[ Y_j^2 | \mathcal{F}_{ j-1 } ]  )  < \frac14-   \sum_{j=1}^k \mathbb{E}[ Y_j^2 | \mathcal{F}_{ j-1 } ]  \Bigg) \nonumber \\
 &\leq& \mathbb{P}\Bigg(   \sum_{j=1}^k    (  Y_j^2 -\mathbb{E}[ Y_j^2 | \mathcal{F}_{ j-1 } ]  )  < \frac14-  (1-\psi(m))  \sum_{j=1}^k \mathbb{E}Y_j^2   \Bigg) \nonumber \\
&\leq& \mathbb{P}\Bigg( \sum_{j=1}^k     Y_j^2-\mathbb{E}[ Y_j^2 | \mathcal{F}_{ j-1 } ]  < -\frac12n  \Bigg).  \label{thn32s}
\end{eqnarray}
Notice  that $$\mathbb{E}[ Y_j^2 | \mathcal{F}_{ j-1 } ] -Y_j^2\leq  \mathbb{E}[ Y_j^2| \mathcal{F}_{ j-1 } ] \leq (1-\psi(m)) \mathbb{E} Y_j^2 \asymp - m.$$
By an argument similar to the proof of (\ref{ine11.1}), we have
$$ \mathbb{E}[  |Y_j^2-\mathbb{E}[ Y_j^2 | \mathcal{F}_{ j-1 } ]|^{(2+\rho)/2}   | \mathcal{F}_{ j-1 } ] \preceq  m^{1+\rho/2}  .$$
Applying Lemma \ref{th24} to  $(\mathbb{E}[ Y_j^2 | \mathcal{F}_{ j-1 } ]-Y_j^2)_{1\leq j \leq k} $  with $\beta=(2+  \rho)/2, x= n/2$ and $v^\beta= km^\beta ,$ from (\ref{thn32s}),  we get
 \begin{eqnarray}
  \mathbb{P}\Bigg(   \sum_{j=1}^k  Y_j^2 < \frac1 4 \Bigg)    & \leq &  \mathbb{P}\Bigg( \sum_{j=1}^k  \mathbb{E}[ Y_j^2 | \mathcal{F}_{ j-1 } ] -  Y_j^2 > \frac12n  \Bigg)  \nonumber \\
  & \leq &  \exp\Big\{ - C(\rho)\, n^{1-\alpha} \Big\},  \label{thks30}
\end{eqnarray}
where $C(\rho)$ is a positive constant.
Notice that
$e^x+ z e^y \leq e^{ z+ x+y}$ for $x, y, z \geq0$ and $z\leq y.$
We obtain  the following upper bound for the relative error of normal approximation:  for all $0\leq x =o(n^{(1-\alpha)/2})  ,$
\begin{eqnarray*}
&&\frac{\mathbb{P}(W_n^o  \geq x )}{1-\Phi \left( x\right)} \leq \frac{\mathbb{P}(W_n^o  \geq x,   \sum_{j=1}^k  Y_j^2 \geq 1/4 )  +  \mathbb{P}( W_n^o  \geq x,    \sum_{j=1}^k  Y_j^2 < 1/4 ) }{1-\Phi \left( x\right)}   \\
&& \leq  \frac{\mathbb{P}( W_n \geq x-C\varepsilon_n, \sum_{j=1}^k  Y_j^2 \geq 1/4 ) + \mathbb{P}(  \sum_{j=1}^k  Y_j^2 < 1/4 ) }{1-\Phi \left( x\right)}\\
&&  \leq  \frac{\mathbb{P}( W_n \geq x-C\varepsilon_n  )  }{1-\Phi \left( x-C\varepsilon_n\right)}   \frac{1-\Phi \left( x-C\varepsilon_n\right) }{1-\Phi \left( x\right)} +   \frac{ \mathbb{P}( \sum_{j=1}^k  Y_j^2 < 1/4 ) }{1-\Phi \left( x\right)}.
\end{eqnarray*}
By (\ref{ineq413}) and (\ref{thks30}), we have for all $0\leq x =o(n^{(1-\alpha)/2})  ,$
\begin{eqnarray*}
&&\frac{\mathbb{P}(W_n^o  \geq x )}{1-\Phi \left( x\right)}  \\
&& \leq \exp\Bigg\{ c_{ \rho} \,  \Bigg(  \frac{  x ^{2+\rho }}{n^{(1-\alpha)\rho / 2  }  }  + x^2 \delta_n^2 + (1+ x) \Big( \frac{1}{n^{(1-\alpha)\rho(2-\rho)/8}(1+x^{\rho(2+\rho)/4})  } +  \delta_n +\varepsilon_n\Big)\Bigg) \Bigg\} \\
&& \ \ \ \ +   \frac{1}{1-\Phi \left( x\right)}   \exp\Big\{ - C(\rho)\, n^{1-\alpha} \Big\} \\
&& \leq \exp\Bigg\{ c'_{ \rho} \,  \Bigg(  \frac{  x ^{2+\rho }}{n^{(1-\alpha)\rho / 2  }  }  + x^2 \delta_n^2 + (1+ x) \Big( \frac{1}{n^{(1-\alpha)\rho(2-\rho)/8}(1+x^{\rho(2+\rho)/4})  } +  \gamma_n \Big)\Bigg) \Bigg\} ,
\end{eqnarray*}
where
\begin{eqnarray*}
  \gamma_n
    =   \delta_n+\varepsilon_n
    \asymp   k^{1/2}  \psi^{1/2}(m ) + n\psi(m )  .
\end{eqnarray*}
Similar, we have the following lower bound for the relative error of normal approximation:  for all $0\leq x =o(n^{(1-\alpha)/2})  ,$
\begin{eqnarray*}
&&\frac{\mathbb{P}(W_n^o  \geq x )}{1-\Phi \left( x\right)}  \nonumber  \\
&&   \geq    \exp\Bigg\{ - c'_{ \rho} \,  \Bigg(  \frac{  x ^{2+\rho }}{n^{(1-\alpha)\rho / 2  }  }  + x^2\delta_n^2 + (1+ x) \Big( \frac{1}{n^{(1-\alpha)\rho(2-\rho)/8}(1+x^{\rho(2+\rho)/4})  } +  \gamma_n \Big)\Bigg) \Bigg\} .
\end{eqnarray*}
Combining the upper and lower bounds of $\frac{\mathbb{P}(W_n^o \geq x )}{1-\Phi \left( x\right)}$ together,
we complete  the proof of Theorem \ref{th3.3}.

\subsection{Proof of Corollary \ref{corollary02}}\label{sec7}
In the proof of Corollary \ref{corollary02}, we will make use of the following well-known inequalities:
\begin{eqnarray}\label{fgsgj1}
\frac{1}{\sqrt{2 \pi}(1+x)} e^{-x^2/2} \leq 1-\Phi ( x ) \leq \frac{1}{\sqrt{ \pi}(1+x)} e^{-x^2/2}, \ \ \ \   x\geq 0.
\end{eqnarray}
First, we show that
\begin{eqnarray}\label{dfgkmsf}
 \limsup_{n\rightarrow \infty}\frac{1}{a_n^2}\ln \mathbb{P}\bigg( \frac{1}{a_n } W_n^o \in B  \bigg) \leq  - \inf_{x \in \overline{B}}\frac{x^2}{2}.
\end{eqnarray}
When $B  =\emptyset,$ the last inequality is obvious. So, we  assume that $B  \neq \emptyset.$
For a given Borel set $B\subset \mathbb{R},$ let $x_0=\inf_{x\in B} |x|.$ Clearly, we have $x_0\geq\inf_{x\in \overline{B}} |x|.$
Therefore, by   Theorem \ref{th3.3},
\begin{eqnarray*}
 \mathbb{P}\bigg(\frac{1}{a_n } W_n^o  \in B \bigg)
 &\leq&  \mathbb{P}\Big(\, \big| W_n  \big|  \geq a_n x_0\Big)\\
 &\leq&  2\Big( 1-\Phi \left( a_nx_0\right)\Big)
  \exp\Bigg\{ c_{ \rho} \,  \Bigg(  \frac{ ( a_nx)^{2+\rho }}{n^{(1-\alpha)\rho / 2  }  }+ \left( a_nx_0\right)^2 \delta_n^2 \\
   &&\ \   + (1+ a_nx) \Big( \frac{1}{n^{(1-\alpha)\rho(2-\rho)/8}(1+(a_nx)^{\rho(2+\rho)/4})  } +  \gamma_n \Big)\Bigg) \Bigg\}.
\end{eqnarray*}
Using   (\ref{fgsgj1}),
we get
\begin{eqnarray*}
\limsup_{n\rightarrow \infty}\frac{1}{a_n^2}\ln \mathbb{P}\bigg(\frac{1}{a_n } W_n^o \in B  \bigg)
 \ \leq \  -\frac{x_0^2}{2} \ \leq \  - \inf_{x \in \overline{B}}\frac{x^2}{2} ,
\end{eqnarray*}
which gives (\ref{dfgkmsf}).

Next, we show that
\begin{eqnarray}\label{dfgk02}
\liminf_{n\rightarrow \infty}\frac{1}{a_n^2}\ln \mathbb{P}\bigg(\frac{1}{a_n }  W_n^o  \in B  \bigg) \geq   - \inf_{x \in B^o}\frac{x^2}{2} .
\end{eqnarray}
When $B^o =\emptyset,$ the last inequality is obvious. So, we  assume that $B^o \neq \emptyset.$
For any given $\varepsilon_1>0,$ there exists an $x_0 \in B^o,$ such that
\begin{eqnarray*}
 0< \frac{x_0^2}{2} \leq   \inf_{x \in B^o}\frac{x^2}{2} +\varepsilon_1.
\end{eqnarray*}
For $x_0 \in B^o$ and all small enough $\varepsilon_2 \in (0, x_0),$ it holds $(x_0-\varepsilon_2, x_0+\varepsilon_2]  \subset B.$
Thus, $x_0\geq\inf_{x\in \overline{B}} x.$ Without loss of generality, we assume that $x_0>0.$
Obviously, we have
\begin{eqnarray*}
\mathbb{P}\bigg(\frac{1}{a_n } W_n^o  \in B  \bigg)   &\geq&   \mathbb{P}\Big(  W_n^o  \in (a_n ( x_0-\varepsilon_2), a_n( x_0+\varepsilon_2)] \Big)\\
&=&   \mathbb{P}\Big(  W_n ^o  \geq  a_n ( x_0-\varepsilon_2)   \Big)-\mathbb{P}\Big( W_n^o  \geq  a_n( x_0+\varepsilon_2) \Big).
\end{eqnarray*}
By Theorem \ref{th3.3}, it is easy to see that $$\lim_{n\rightarrow \infty} \frac{\mathbb{P}\Big( W_n^o  \geq  a_n( x_0+\varepsilon_2) \Big)}{\mathbb{P}\Big(  W_n ^o  \geq  a_n ( x_0-\varepsilon_2)   \Big) } =0 .$$
Then, by (\ref{fgsgj1}), it follows that
\begin{eqnarray*}
\liminf_{n\rightarrow \infty}\frac{1}{a_n^2}\ln \mathbb{P}\bigg(\frac{1}{a_n }  W_n^o \in B  \bigg)  \geq  -  \frac{1}{2}( x_0-\varepsilon_2)^2 . \label{ffhms}
\end{eqnarray*}
Now, letting $\varepsilon_2\rightarrow 0,$  we have
\begin{eqnarray*}
\liminf_{n\rightarrow \infty}\frac{1}{a_n^2}\ln \mathbb{P}\bigg(\frac{1}{a_n }  W_n^o  \in B \bigg) \ \geq\ -  \frac{x_0^2}{2}  \  \geq \   -\inf_{x \in B^o}\frac{x^2}{2} -\varepsilon_1.
\end{eqnarray*}
Because $\varepsilon_1$ can be arbitrarily small, we get (\ref{dfgk02}).
Combining (\ref{dfgkmsf}) and (\ref{dfgk02}) together, we complete  the proof of Corollary \ref{corollary02}.

\end{document}